\documentclass[letterpaper,10pt]{article}
\usepackage{amsthm}
\usepackage{amsmath,amssymb}

\newcommand{\mr}[1]{\mathcal{#1}}

\newtheorem{lm}{Lemma}[section]
\newtheorem{thq}[lm]{Theorem}
\theoremstyle{remark}
\newtheorem{rem}[lm]{Remark}
\newtheorem{df}[lm]{Definition}
\newtheorem{co}[lm]{Corollary}
\newtheorem{que}[lm]{Question}

\DeclareMathOperator{\SL}{SL}
\DeclareMathOperator{\Hom}{Hom}

\newcommand{\la}{\langle}
\newcommand{\ra}{\rangle}

\newcommand{\C}{\mathbb{C}}
\newcommand{\Q}{\mathbb{Q}}
\newcommand{\R}{\mathbb{R}}
\newcommand{\T}{\mathbb{T}}
\newcommand{\Z}{\mathbb{Z}}
\newcommand{\N}{\mathbb{N}}
\newcommand{\F}{\mathbb{F}}

\newcommand{\cl}{\mathcal{L}}

\title{Kazhdan Constants for $\SL_n(\Z)$}
\author{Martin Kassabov\thanks{The author is supported by several NSERC grants}}
\begin{document}
\maketitle

\begin{abstract}
In this article we improve the known Kazhdan constant for
$\SL_n(\Z)$ with respect to the generating set of the elementary
matrices. We prove that the Kazhdan constant is bounded from below
by $[42\sqrt{n}+860]^{-1}$, which gives the exact asymptotic
behavior of the Kazhdan constant, as $n$  goes to infinity, since
$\sqrt{2/n}$ is an upper bound.

We can use this bound to improve the bounds for the spectral gap
of the Cayley graph of $\SL_n(\F_p)$ and for the working time of
the product replacement algorithm for abelian groups.
\end{abstract}

% ------------------------------------------

\begin{section}{Introduction}

%Kazhdan property $T$
%\marginpar{...}

Kazhdan property $T$ plays important role in the representation
theory of arithmetic groups. Since the work of Kazhdan
(see~\cite{Kazhdan}) it is known that any higher rank arithmetic
group has property $T$.

In recent years there have been several connections between
property $T$ and the working time of several algorithms in
componential group theory. In~\cite{LubPak} the authors use the
Kazhdan property $T$ of the group $\SL_n(\Z)$ to prove that the
product replacement algorithm on abelian groups has logarithmic
working time. In order to make these results quantitative one
needs exact values of the Kazhdan constants for certain groups,
most notably for the group $\SL_n(\Z)$, with respect to the
standard generators.

%Importance of explicit Kazhdan constants. \marginpar{...}

Almost all methods for proving that an arithmetic group $\Gamma$
has property $T$ use Kazhdan's result and transfer the question to
Lie groups. These methods can not be made quantitative and they do
not lead to any explicit Kazhdan constants. The first author to
obtain partial results in this direction was M.~Burger ---
in~\cite{Burger}, he found a lower bound for the some constant,
closely related to the  Kazhdan constant for the group
$\SL_3(\Z)$.
%in the class of all finite dimensional representations.
Several years later, Y.~Shalom (see~\cite{Shalom1}) used bounded
generation to construct an explicit  lower bound for Kazhdan
constant for the group $\SL_n(\Z)$. This result can be combined
with the result of O.~Tavgen (see~\cite{Tag}) to obtain similar
bounds for other higher rank arithmetic groups.

The main result of this paper gives the exact asymptotical
behavior of the Kazhdan constant of $\SL_n$ with respect to the
generating set $E_n$ consisting of elementary matrices with %
$\pm 1$ off the diagonal. The methods used in the proof are based
on the ideas in~\cite{Shalom1}.

%Shalom's result

Structure of the paper: in the following section we describe the
main result and an outline of the basic idea of the proof;
section~\ref{sec:applications} we give several applications of the
main theorem; sections~\ref{sec:relp} and~\ref{sec:relpq} are
dedicated to the derivations of explicit Kazhdan constants for the
relative property $T$ of the groups  $\SL_p \ltimes \Z^{p}$ and %
$(\SL_p \times \SL_q) \ltimes \Z^{pq}$; section~\ref{sec:vec}
describes vector systems in $\Z^k$ and generalized elementary
operations, which are used in the proof of the stronger version
bound generation property of the group $\SL_n(\Z)$ with respect to
the set $E_n$; section~\ref{sec:proof} concludes the proof of
Theorem A. The last section is dedicated to some possible
extensions of the main theorem.
%\marginpar{...}

Acknowledgements: The author wishes to thanks Roman Muchnik, Tal
Poznansky and Misha Ershov for the useful discussions and helpful
suggestions during the preparation of the manuscript. I
thank Igor Pak for suggesting several applications of the main
result. I wish to express my gratitude to Alex Lubotzky, Yehuda Shalom
and my adviser Efim Zelmanov for introducing me to this subject.
I also wish to thank to Clay Mathematics Institute for the
financial support during the preparation of this paper.
%
%Roman for measures on $\T^p$;
%Tal for primes subgroups and the bound $2k+1$ in ..;
%Misha for K-theory and proving me with some useful references
%Shalom for suggestions
%\marginpar{Clay Institute}

\end{section}

\begin{section}{Main Result}
\label{sec:main}

Let us recall the definition of Kazhdan property $T$:
\begin{df}
A topological group $G$, generated by a compact set $Q$, is said
to have Kazhdan property $T$, if there exists a constant
$\epsilon$, such that any (continuous) unitary representation
$(\pi,\mr{H})$ of the group $G$, which contains a unit vector $v$
such that $||\pi(g)v-v|| \leq \epsilon$ for any $g\in Q$, contains
a $G$ invariant vector. The maximal $\epsilon$ with this property
is called the Kazhdan constant of $G$ with respect to $Q$ and is
denoted by $\mr{K}(G,Q)$.
\end{df}

In~\cite{Kazhdan}, Kazhdan proved that any higher rank Lie group
$G$ and any lattice $\Gamma$ in such a group has property $T$,
with out giving any values for the Kazhdan constants.\footnote{
This result does not hold for some rank $1$ groups. For example,
the groups $\SL_2(\R)$ and $\SL_2(\Z)$ does not have Kazhdan
property $T$.
}
In
particular from his work follows that $\SL_n(\Z)$ has property $T$ if
$n\geq 3$.

Let $E_n$ be the set of all elementary matrices with $\pm 1$ off
the diagonal. It is known that the set $E_n$ generates the group
$SL_n(\Z)$ and it is natural to look for the %exact
value of the
Kazhdan constant $\mr{K}(\SL_n,E_n)$. The main result in this
paper is the following lower bound for the Kazhdan constant of the
group $\SL_n(\Z)$ with respect to the set $E_n$.

%This result improves several bound for the working times of
%important algorithms in computational group theory.

{\renewcommand{\thelm}{A} \addtocounter{lm}{-1}
\begin{thq}
%\marginpar{...} %
\label{main} The Kazhdan constant for $\SL_n(\Z)$ with respect to
the set $E_n$ of elementary matrices with $\pm 1$ off the diagonal
is %(the bound can be slightly improved)
$$
\mr{K}(\SL_n(\Z),E_n) \geq (64\sqrt{n}+2850)^{-1}.
$$
\end{thq}
}

Using the same methods but carefully tracking all constants,
allows us to obtain a slightly better result:

{\renewcommand{\thelm}{A'} \addtocounter{lm}{-1}
\begin{thq}
%\marginpar{...} %
\label{main:strong} The Kazhdan constant for $\SL_n(\Z)$ with
respect to the set $E_n$ of elementary matrices with $\pm 1$ off
the diagonal is
$$
\mr{K}(\SL_n(\Z),E_n) \geq (42\sqrt{n}+860)^{-1}.
$$
\end{thq}
}

If we consider the group $\SL_n(\F_p)$ we can improve the bound
even further:
{\renewcommand{\thelm}{A''} \addtocounter{lm}{-1}
\begin{thq}
%\marginpar{...} %
\label{main:fp} The Kazhdan constant for $\SL_n(\F_p)$ with
respect to the set $E_n$ of elementary matrices
with $\pm 1$ off the diagonal is %(the bound can be slightly improved)
$$
\mr{K}(\SL_n(\F_p),E_n) \geq (31\sqrt{n}+700)^{-1}.
$$
\end{thq}
}

Theorem A can be generalized to linear group over number fields:
{\renewcommand{\thelm}{B} \addtocounter{lm}{-1}
\begin{thq}
%\marginpar{...} %
\label{main:nf}
Let $\mr{O}$ be the ring of integers in a number field $\mathbb{K}$,
%with discriminant $\Delta$,
which is generated as a ring by $1$ and the elements $\alpha_i$
for $i=1,\dots,s$. The Kazhdan constant for $\SL_n(\mr{O})$ with
respect to the set $E_n(\mr{O}$) of  elementary matrices with %
$\pm 1$ and $\pm \alpha_i$ off the diagonal is
$$
\mr{K}(\SL_n(\mr{O}),E_n(\mr{O})) \geq
   [50\sqrt{n} + (B + C\Delta)6^m]^{-1}.
$$
where $B$ and $C$ are universal constants and $\Delta$ is the
number of different prime divisors of the discriminant of
$\mathbb{K}:\Q$.\footnote{From the Generalized Reimann Hypothesis
it follows that the bound for the Kazhdan constant does not depend
on the discriminant of $\mathbb{K}:\Q$, see~\cite{Carter2} for
details. } %\marginpar{footnote}
\end{thq}
}

The lower bounds of the Kazhdan constant obtained in Theorems~A
and B are asymptotically exact. Let us consider the natural
representation of $\SL_n(\Z)$ into $\cl^2(\Z^n)$.
% which does not have invariant vectors.
There exists a unit vector %
$v\in \cl^2(\Z^n)$, which is moved by any element in $E_n$ by a
distance of $\sqrt{2/n}$.\footnote{This representation contains
one dimensional space of invariant vectors, however the vector $v$
lies in the orthogonal compliment of this invariant subspaces and
we cans restrict every thing to that subspace.}
 This shows that
$$
\mr{K}(\SL_n(\Z),E_n) \leq \sqrt{2/n}.
$$
This upper bound for the Kazhdan constant for $\SL_n(\Z)$ was
found by A.~Zuk (unpublished) and  can be found
in~\cite{Shalom1}.
%$ is an upper bound for the Kazhdan constant.

The basic idea of the proof is similar to the one in
\cite{Shalom1}. There, Y.~Shalom used the relative property $T$
and bounded generation to prove that $(33n^2-11n + 1152)^{-1}$ is
a lower bound for the Kazhdan constant. The main steps in his
proof are the following:

The group $\SL_2(\Z) \ltimes \Z^2$ has relative property $T$, with
respect to the normal subgroup $\Z^2$ with Kazhdan constant at
least $1/10$ (using the
elementary matrices in $\SL_2$ and the basis vectors of $\Z^2$ as
generating set). This bound gives that for any unitary
representation $(\pi,\mr{H})$ if any elementary matrix %
(with $\pm 1$ off the diagonal) moves a fixed unit vector %
$v \in \mr{H}$ by less than $\epsilon$, then any elementary matrix
(with any integer off the diagonal) moves the same vector $v$ by at
most $20\epsilon$.

Bounded generation of the group $\SL_n(\Z)$ with respect to the
elementary matrices, proved by Carter and Keller~\cite{Carter2},
gives that any element $g\in \SL_n(\Z)$ can be written as product
of at most (approximately) $3n^2/2$ elementary matrices. This,
together with the previous step, shows that any element %
$g\in \SL_n(\Z)$ moves the vector $v$ by at most $30n^2\epsilon$.

Finally we use the observation that if a unit vector is moved by any
element of a group by a distance less than $1$ then the
representation has an invariant vector. This leads to a lower
bound for the Kazhdan constant of $SL_n(\Z)$ of the type
$O(n^{-2})$.

Using  this idea is not possible to obtain an asymptotically better
bound, because a generic element in $\SL_n(\Z)$ can not be written
as a product of less than $n^2$ elementary matrices. In our
proof, instead of working with the group $\Z^2$, we work with larger
abelian subgroups $H_i$ of $SL_n$.
% (not only with $\Z$).
Instead of $\SL_2 \ltimes \Z^2$, we use the group
$N_i\ltimes H_i$, where $N_i$ is semi-simple and $N_i\ltimes H_i$
is maximal parabolic. This group has a relative property $T$ with
respect to the subgroup $H_i$ and in
section~\ref{sec:relpq} we find a lower
bound for the relative Kazhdan constant.%
\footnote{In almost all cases, the group $N_i\ltimes H_i$ actually
has property $T$. Using the same methods a bound the Kazhdan
constant for this group can be computed, and it is of the similar
to the relative Kazhdan constant.} %
Using this constant it can be shown that if any elementary matrix
move $v$ by at least $\epsilon$, then any element $g$ lying in
some $H_i$, moves $v$ by at most $2k(n) \epsilon$. We can
obtain a this result using the relative property $T$ of
$\SL_2(\Z)\ltimes \Z^2$ which can be embedded in many different ways in
$N_i \ltimes H_i$ and will obtain an upper bound for $k(n)$ of type
$O(n^2)$.
However, using the whole groups $N_i \ltimes H_i$ allows us to improve the
bound for $k(n)$ to $O(\sqrt{n})$, which allows us to obtain a
better lower bound for
the Kazhdan constant of $\SL_n(\Z)$.

%The function $k(n)$ is not bounded, but $k(n)$ increase slower
%than $n^2$ which is the dimension of a generic $H_i$.

The main result in section~\ref{sec:vec} shows that using 5
multiplications by elements lying in conjugates of some $H_i$,
every element in $\SL_n(\Z)$ can be transformed to an element in
$\SL_{\lambda n}(\Z)$ (embedded in the upper left corner of
$\SL_n(\Z)$), where $\lambda = 2/3$. This, together with the result
by Carter and Keller~\cite{Carter,Carter2}, saying that any matrix
in $\SL_3(\Z)$ can be written as a product of at most $60$ elementary
matrices, gives that any matrix in $\SL_n(\Z)$ can be written as a
product of at most $s(n)=60 + 13\ln n$ matrices lying in
some $H_i$.%
\footnote{Using a result of L. Vaserstain (see~\cite{Vas}) it can
be shown that any element in $\SL_n(\Z)$ can be written as a
product of a fixed number of matrices in some
$H_i$. The best bound for the number of matrices need to write any
element in $SL_n$ as such product is around 100 and leads to
slightly worse bound for the Kazhdan constant of $\SL_n(\Z)$ than
the one obtained in Theorem~\ref{main}.}
This implies that any element in %
$\SL_n$ moves the vector $v$ by at most  $h(n) \epsilon$, where
$h(n) = 2s(n)k(n)$. From this result, we obtain that the Kazhdan
constant for $\SL_n(\Z)$ with respect to the set of elementary
matrices is at least $1/h(n)$. This argument gives a bound of type
$(\sqrt{n}\ln n)^{-1}$ for the Kazhdan constant, but a more detailed
consideration allow us to improve this bound to $n^{-1/2}$.

\end{section}

\begin{section}{Applications}
\label{sec:applications}

%Application of our result. \marginpar{...}

%Using the results obtained by Pak and Zuk in~\cite{PakZuk}, we can obtain the
%following corollaries of our result
%\begin{list}{}{}

The value of the Kazhdan constant is related to several other
constants, like the spectral gap of the Laplacian and the
mixing time of random walks on finite Cayley graphs.

\begin{subsection}{Spectral gap of Cayley graph of $\SL_n(\Z)$}
Applying the result by I. Pak and A. Zuk
from~\cite{PakZuk}, we have that the spectral
gap is
% \marginpar{ask Pak if we can improve this to 2}
$$
 \beta(\SL_n(\Z)) \geq \mr{K}(\SL_n(\Z))^2/4
 \geq \frac{1}{4(42\sqrt{n} + 860)^2} = O(1/n).
$$
This bound is better than the previously known one which was
$O(n^{-4})$. The argument from section~\ref{sec:main} give that
$$
\beta(\SL_n(\Z)) \leq 1/n.
$$
\end{subsection}

\begin{subsection}{Spectral gap of Cayley graph of $\SL_n(\F_p)$}
Since $\SL_n(\F_p)$ is a factor group of $\SL_n(\Z)$, the spectral
gap of $\SL_n(\F_p)$ is bigger than the one for $\SL_n(\Z)$. We
can obtain a slightly better estimate, using the fact that $\F_p$ is
a field.
\iffalse
and we can transform any vector system (over $\F_p$) to
the standard one using at most $3$ generalized elementary
transformations (if $n\geq 2k+1$). This gives a
\fi
The better bound for
the Kazhdan constant
$$
\mr{K}(\SL_n(\F_p)) \geq (31\sqrt{n}+700)^{-1},
$$
from Theorem~A', yields the bound for the spectral gap
$$
\frac{1}{n} \geq  \beta(\SL_n(\F_p)) \geq
     \frac{1}{4(31\sqrt{n} + 700)^2} = O(1/n).
$$
\end{subsection}

\begin{subsection}{Mixing time of random walks on $\SL_n(\F_p)$}
There is a connection between the spectral gap of the finite
Cayley graph and the mixing time of the random walk on the same
graph. Applying that to the group $G= \SL_n(\F_p)$ gives
$$
mix \lesssim \beta \log |G|  = O(n^3 \log p),
$$
which is better than the previous known bounds of $n^6 \log p$
(see~\cite{PakZuk}) and $n^4 \log^3 p$ (see~\cite{DiaCos}).
\end{subsection}

\iffalse
\begin{subsection}{Diameter of the Cayley graph for $\SL_n(\F_p)$}

\end{subsection}
\fi

\begin{subsection}{Mixing time of the product
replacement algorithm for abelian groups}
In~\cite{LubPak},
A. Lubotzky and I. Pak showed a connection between the
working time of the Product Replacement Algorithm on
$n$ generated abelian groups and the Kazhdan constant for $SL_n(\Z)$,
in particular they proved the following upper bound for the
mixing time:
\iffalse
the following bound of the mixing time of the
product replacement algorithm for $n$ generated abelian group is
proved
\fi
$$
mix \lesssim n  \mr{K}(\SL_n(\Z),E_n)^{-2} \log |\Gamma|.
$$
Using the bound for the Kazhdan constant from
Theorem~\ref{main} gives a bound of $n^2 \log |\Gamma|$
\end{subsection}

\begin{subsection}{Relaxation time for particle systems}
Particle systems was studied by P.~Diaconis and L.~Saloff-Coste
(the original problem was proposed by D. Aldous). In~\cite{DiaCos}
they proved that the relaxation time of a particle system is
bounded by $n^2\log n$ and made a conjecture that the bound is of
type $n\log n$. Vies the particle system as a random walk on
$\Z_2^n$ we can use the bound for the Kazhdan constant of
$\SL_n(\Z)$, which allows us to confirm this conjecture.
\end{subsection}

%\end{list}
\end{section}

\begin{section}{Relative Kazhdan constant for
$\SL_p(\Z)\ltimes \Z^{p}$}
\label{sec:relp}

In this section, we estimate the relative Kazhdan constant for the
group $\SL_p(\Z) \ltimes \Z^{p}$ with respect to the group
$\Z^{p}$, considering the set of elementary matrices in $\SL_p$
together with the basis vectors of $\Z^{p}$ as a generating set.
The idea of the proof of Theorem~\ref{relp} is based to the one used
by Burger in~\cite{Burger}%
\footnote{In his paper~\cite{Burger}, Burger attributes the idea
of this proof to Furstenberg. } %
for estimating the relative Kazhdan constant for %
$\SL_2(\Z) \ltimes \Z^2$. A very detailed explanation of this
proof can be found in~\cite{Shalom1}.

\begin{thq}
\label{relp}
Let $F$  denote the set of elementary matrices in
$\SL_p(\Z)$, and
$G$ denote the set of the $p$ standard basis elements of $\Z^{p}$.
Let $(\pi,\mr{H})$ be
a unitary representation of $\SL_p(\Z)\ltimes \Z^{p}$,
containing a vector $v$ which is $(F\cup G, 1/l(p))$
invariant, where %\marginpar{check it!}
$$
l(p) = \sqrt{p+25}+3.
$$
Then $\mr{H}$ contains $\Z^{p}$ invariant vector, provided that
$p\geq 2$.
\end{thq}
\begin{rem}
Let us consider the standard the unitary representation of
$\SL_{p+1}$ on $\cl^2(\Z^{p+1})$. The group $\SL_p\ltimes \Z^{p}$
is isomorphic to a maximal parabolic in $\SL_{p+1}$, thus we have
a representation of that group in $\cl^2(\Z^{p+1})$. The
representation decomposes as a sum of two representations, one is
trivial and isomorphic to $\cl^2(\Z)$, the other is a
representation on $\cl^2(\Z^{p+1} \setminus \Z)$, without $\Z^{p}$
invariant vectors. Using that representation and a suitable vector
$v$, it can be shown that the Kazhdan constant is at most
$\sqrt{2/p}$. This shows that $1/\tilde{l}(n)$ is not a relative
Kazhdan constant, where $\tilde{l}(p) =\sqrt{p/2}$.
\end{rem}
\iffalse
\begin{rem}
If we simply generalized the proof by Burger, we will get a bound
of the type $1/p$. This is true because if we consider the action
of $\SL_p$ on the projective space $P\R^{p-1}$, there exists
probability measure $\mu$ such that for every Borel set $B$ we
have $|\mu(gB) - \mu(B)| \leq 1/{p}$. The constant in
theorem~\ref{relp} is better because we use
lemma~\ref{invariance}, which is stronger than the similar results
used in~\cite{Burger} and in~\cite{Shalom1}.
\end{rem}
\fi

\begin{proof}
Assume that $v\in\mr{H}$ is $(F\cup G, \epsilon)$ invariant
unit vector, and that the Hilbert space $\mr{H}$ does not contain
$\Z^{p}$ invariant vector.
%Let us restrict the representation $\pi$ to the $\Z^{pq}$ and
Let $P$ be the % corresponding
projection valued measure on $\widehat{\Z^{p}} = \T^{p}$, coming
form the restriction of the representation $\pi$ to $\Z^{p}$, and
let $\mu_v$ be the measure on $\T^{p}$, defined by %
$\mu_v(B)=\la P(B)v,v \ra$. The probability measure
$\mu_v$ is supported on
$\T^{p}\setminus \{0\}$, because by assumption $\mr{H}$ does not
contain an $\Z^{p}$ invariant vector and by construction
$P(\{0\})$ is the
projection onto the space of $\Z^p$ invariant vectors.
%which implies that $P(\{0\}) = 0$.

For an element $x \in \T^{p}$ we will write $x=(x_{1},\dots,x_{p})$,
where all $x_{i}$ are in $\R/\Z$, which
we identify with the interval $(-1/2,1/2]$.

\begin{lm}
Let $K_{i} = \{x \mid 1/4 > |x_{i}| \}$, then %
$\mu_v(K_{i}) \geq 1-\epsilon^2/2$.
\end{lm}
\begin{proof}
By the definition of the measure $\mu_v$, we have
$$
|| \pi(g_{i})v-v||^2 =
\int_{\T^{p}} | e^{2\pi i x_{j}} -1 | ^2\, d\mu_v \leq \epsilon^2,
$$
where $g_{i}$ form the standard basis of $\Z^{p}$. Now using the
fact that $| e^{2\pi i x_{i}} -1 |^2 \geq 2$ for $1/2 \geq |x_{i}|
\geq 1/4$, the above inequality implies that $\mu_v(|x_{i}| \geq
1/4) \leq \epsilon^2/2$.
\end{proof}

\begin{lm}
\label{invariance} For every Borel set $B \subset \T^{p}$ and
every elementary matrix $g \in F$, we have that
$$
|\mu_v(gB) - \mu_v(B)| \leq 2\epsilon\sqrt{\mu_v(B)}+ \epsilon^2.
$$
The action of $\SL_p(\Z)$ on $\T^{p}$ is the standard one coming
from the action on $\R^{p}$, via the isomorphism
$\T^{p}=\R^{p}/\Z^{p}$.
\end{lm}
\begin{rem}
Similar lemma was used in~\cite{Burger} and in~\cite{Shalom1}, but
the upper bound for $ |\mu_v(gB) - \mu_v(B)|$ was $2\epsilon$. If
we use that lemma we could only obtain $1/p$ as a bound for the
relative Kazhdan constant for $\SL_p \ltimes \Z^p$, which will
give a Kazhdan constant for $\SL_n(\Z)$ of the form $O(1/n)$.
\end{rem}
\begin{proof}
Using the properties of the projection valued measure $P$, we have
$$
\begin{array}{r@{\,\,}l}
 |\mu_v(gB) - \mu_v(B)| = &
       | \la\pi(g^{-1})P(B)\pi(g)v,v\ra - \la P(B)v,v\ra| \leq \\
 \leq &
       |\la\pi(g^{-1})P(B)(\pi(g)v-v),v\ra | +
       |\la P(B)v,(\pi(g)v-v\ra)| = \\
 = &
       2|\la\pi(g)v-v,P(B)v\ra| +
       \la P(B)(\pi(g)v-v),\pi(g)v-v\ra \leq  \\
 \leq &
       2\epsilon \sqrt{ \mu_v(B)} + \epsilon^2,
\end{array}
$$
where the final inequality follows from the facts
that $v$ is $(F,\epsilon)$ invariant
vector and $||P(B)v||^2 = \mu_v(B)$.
\end{proof}

\begin{lm}
\label{R2} Let $\mu$ be a finitely additive measure on $\T^2$ such
that:
\begin{itemize}
\item %
$\mu(|x| \geq 1/4) \leq \epsilon^2/2$ and %
$\mu(|y| \geq 1/4) \leq \epsilon^2/2$, %
\item %
$|\mu(gB) - \mu(B)| \leq 2\epsilon \sqrt{\mu(B)} + \epsilon^2$ %
for any Borel set $B$ and any elementary matrix $g\in \SL_2(\Z)$.
\end{itemize}
Then
%if $\epsilon < 1/2$
we have
$$
\mu(\T^2\setminus \{(0,0)\}) \leq (2+\sqrt{10})^2\epsilon^2
\quad \mbox{and} \quad
\mu(x\not=0, y=0) \leq (1+\sqrt{3})^2\epsilon^2.
$$
\end{lm}
\begin{proof}
Let us define the Borel subsets $A_i$ and $A_i'$ of $\T^2$ using
the picture:

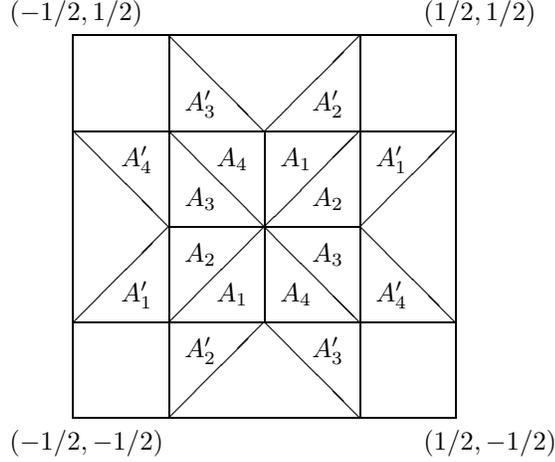
\setlength{\unitlength}{3947sp}%
\begin{picture}(3000,2900)(200,-2350)
\thinlines
\put(1800,-2100){\framebox(2400,2400){}}
\put(1800,-1500){\line(1,0){2400}}
\put(1800, -300){\line(1,0){2400}}
\put(2400, -900){\line(1,0){1200}}
\put(2400,-2100){\line(0,1){2400}}
\put(3600,-2100){\line(0,1){2400}}
\put(3000,-1500){\line(0,1){1200}}
\put(2400,-1500){\line(1, 1){1200}}
\put(2400, -300){\line(1,-1){1200}}
\put(1800,-1500){\line(1, 1){600}}
\put(1800, -300){\line(1,-1){600}}
\put(3600, -900){\line(1, 1){600}}
\put(3600, -900){\line(1,-1){600}}
\put(2400,-2100){\line(1, 1){600}}
\put(2400,  300){\line(1,-1){600}}
\put(3000,-1500){\line(1,-1){600}}
\put(3000, -300){\line(1,1){600}}
\put(3100,-520){$A_1$}
\put(3300,-770){$A_2$}
\put(3300,-1120){$A_3$}
\put(3100,-1370){$A_4$}
\put(2700,-1370){$A_1$}
\put(2500,-1120){$A_2$}
\put(2500,-770){$A_3$}
\put(2700,-520){$A_4$}
\put(3700,-520){$A_1'$}
\put(3300,-170){$A_2'$}
\put(3300,-1720){$A_3'$}
\put(3700,-1370){$A_4'$}
\put(2100,-1370){$A_1'$}
\put(2500,-1720){$A_2'$}
\put(2500,-170){$A_3'$}
\put(2100,-520){$A_4'$}
\put(1400,-2300){$(-1/2,-1/2)$}
\put(1400, 400){$(-1/2,1/2)$}
\put(4000,-2300){$(1/2,-1/2)$}
\put(4000,400){$(1/2,1/2)$}
\end{picture}

\noindent
%\marginpar{is this enough?}
Each set $A_i$ or $A_i'$ consists of the interiors of two
triangles and part of their boundary (not including the vertices).
The sets $A_i$ do not contain the side which is part of the small
square, they also do not contain their clockwise boundary
but contain the counter-clockwise one. Each set $A_i'$ includes
only the part of its boundary which lies on the small square.

From the picture it can be seen that the elementary matrices
$g_{ij}^{\pm} = I \pm e_{ij}\in F$,
act on the sets $A_i$ as follows:
$$
\begin{array}{ll}
g_{12}^+(A_3 \cup A_4')  = A_3 \cup A_4 \quad &
g_{21}^+(A_3'\cup A_4 )  = A_3 \cup A_4 \\
g_{12}^-(A_1'\cup A_2 )  = A_1 \cup A_2  \quad &
g_{21}^-(A_1 \cup A_2')  = A_1 \cup A_2.
\end{array}
$$
Using the properties of the measure $\mu$ the above equalities
imply the inequalities:
$$
\begin{array}{l}
\mu(A_1) + \mu(A_2) \leq
     \mu(A_1') + \mu(A_2) + \epsilon^2 +
     2\epsilon\sqrt{\mu(A_1') + \mu(A_2 )} \\
\mu(A_1) + \mu(A_2) \leq
     \mu(A_1) + \mu(A_2') + \epsilon^2 +
     2\epsilon\sqrt{\mu(A_1 ) + \mu(A_2')} \\
\mu(A_3) + \mu(A_4) \leq
     \mu(A_3') + \mu(A_4) + \epsilon^2 +
     2\epsilon\sqrt{\mu(A_3') + \mu(A_4 )} \\
\mu(A_3) + \mu(A_4) \leq
     \mu(A_3) + \mu(A_4') + \epsilon^2 +
     2\epsilon\sqrt{\mu(A_3 ) + \mu(A_4')}. \\
\end{array}
$$
Adding these inequalities and noticing that
$$
\begin{array}{l}
\mu(A_1') + \mu(A_4') \leq
     \mu(\{ |x| \geq 1/4\}) \leq
     \epsilon^2/2 \,\,\, \mbox{ and}\\
\mu(A_2') + \mu(A_3') \leq
     \mu(\{ |y| \geq 1/4\}) \leq
     \epsilon^2/2
\end{array}
$$
we obtain
$$
\sum_i \mu(A_i) \leq
    4\epsilon^2 + \sum_i \mu(A_i') +
        2\epsilon \sqrt{4\left(\sum \mu(A_i) +
    \sum \mu(A_i')\right)} \leq
$$
$$
    \leq 5\epsilon^2 + 4\epsilon \sqrt{\sum \mu(A_i) + \epsilon^2}.
$$
%in order to obtain the firs inequality,
Here we have used
that any positive numbers $a_i$, satisfy the inequality
$$
\sum_{i=1}^k \sqrt{a_i} \leq \sqrt {k \sum_{i=1}^k a_i}.
$$

After substituting $c =\sqrt{\sum \mu(A_i) + \epsilon^2}$
and solving the resulting quadratic inequality we obtain
$\sum\mu(A_i) \leq (13 +4 \sqrt{10}) \epsilon^2$.

Also from the system of inequalities, taking the inequality for
$\mu(A_i)$, where the index $i$ is such that $\mu(A_i)$ is maximal,
we have:
$$
\max_i \mu(A_i) \leq %\rule[0pt]{15pt}{0pt}
    \epsilon^2 + \max_i \mu(A_i') +
         2\epsilon \sqrt{\max_i \mu(A_i) + \max_i \mu(A_i')} \leq
$$
$$
    \leq 3\epsilon^2/2 +
        2\epsilon \sqrt{ \max_i \mu(A_i) + \epsilon^2/2},
$$
which yields
$\max_i \mu(A_i) \leq (7/2 + 2\sqrt{3}) \epsilon^2$.
Finally we can use that
$$
\mu(\T^2\setminus \{(0,0)\}) \leq \sum \mu(A_i) + \mu(\{ |x| \geq
1/4\}) + \mu(\{ |y| \geq 1/4\}) \leq
$$
$$
\leq (14 + 4 \sqrt{10}) \epsilon^2 = (2+\sqrt{10})^2\epsilon^2,
$$
and
$$
\mu(x\not=0,y=0) \leq \max_i \mu(A_i) + \mu(\{ |x| \geq 1/4\}) \leq
%\frac{4\epsilon^2-2\epsilon^3}{2(1-4\epsilon)} =
(4 + 2 \sqrt{3})\epsilon^2 = (1+\sqrt{3})^2 \epsilon^2,
$$
which completes the proof of the lemma.
% |\/|
%-+--+-
%\|\/|/
%/|/\|\
%-+--+-
% |/\|
\end{proof}
\begin{lm}
\label{Rp}
Let $\mu$ be a finitely additive measure on $\T^p$,
which satisfies conditions from the previous lemma (with $\SL_p$
replacing $\SL_2$) %. Then, if $\epsilon \leq 1/4$,
$$
\mu(\T^p\setminus \{(0,\dots,0)\}) \leq
(\sqrt{p+25} + 3)^2 \epsilon^2.
$$
\end{lm}
%Change the proof and make it similar to the one of the next lemma.
\begin{proof}
For a point $y\in \T^p$ we write $y=(y_{1}\dots,y_{p})$, where
$y_i \in (1/2,1/2]$. Let us define the Borel sets
$$
\begin{array}{l}
B_i = \{ y \mid y_k =0 \mbox{ for } k\leq i\}, \mbox{ and}\\
C_i = \{ y \mid y_1 = y_i \not =0, y_k = 0 \mbox{ for } 1<k<i\}.
\end{array}
$$
The elementary matrix $g_{1i} \in \SL_p$ sends $B_{i-1} \setminus
B_i$ into $C_i$ for any $i\geq 3$. Therefore, we have
$$
\mu(B_{i-1} \setminus B_i ) \leq
\mu(C_i) + \epsilon^2 + 2\epsilon\sqrt{\mu(C_i)}.
$$
Let us notice that the sets $C_i$, for $i=2,\dots, p$, are
disjoint and their union lies in the set %
$C=\{ y | y_1 \not =0, y_2 = 0 \}$. Therefore by
adding these inequalities we have
$$
\begin{array}{r@{\,\,}l}
\mu(B_2 \setminus B_p ) = & \displaystyle
\sum_{i=3}^p \mu(B_{i-1}\setminus B_i) \leq \\
\leq & \displaystyle
\sum_{i=3}^p [\mu(C_i) + \epsilon^2 + 2\epsilon\sqrt{\mu(C_i)}] \leq \\
\leq & \displaystyle
\mu(\cup_i C_i) + (p-2)\epsilon^2 +
    2\epsilon\sqrt{(p-2)\mu(\cup_i C_i)}\leq \\
\leq & \displaystyle \rule[15pt]{0pt}{0pt}
 \mu(C) + (p-2)\epsilon^2 + 2\epsilon\sqrt{(p-2)\mu(C)}.
\end{array}
$$
Using the projection $\T^p \to \T^2$ given by taking at the
first two coordinates, we can project the measure $\mu$ to a
measure $\tilde \mu$ on $\T^2$. Applying the previous lemma
to the measure $\tilde \mu$ we have:
$\displaystyle \mu(\T^p \setminus B_2) \leq
    (2+\sqrt{10})^2\epsilon^2$ and
$\displaystyle \mu(C) \leq
    (1+\sqrt{3})^2\epsilon^2$,
therefore
$$
\begin{array}{r@{\,\,}l}
\mu(\T^p\setminus & \{(0,\dots,0)\}) =
    \mu(\T^p \setminus B_2) + \mu(B_2 \setminus B_p) \leq \\
\leq & \displaystyle
    (2+\sqrt{10})^2\epsilon^2 + (1+\sqrt{3})^2\epsilon^2 +
    (p-2)\epsilon^2 +
    2(1+\sqrt{3})\sqrt{p-2} \epsilon^2 =\\
= & \displaystyle
    \left(p+16 + 4\sqrt{10} + 2\sqrt{3}
        +2(1+\sqrt{3})\sqrt{p-2}\right)\epsilon^2 \leq \\
\leq & (p+ 6\sqrt{p} + 33)\epsilon^2
\leq (\sqrt{p+25} + 3)^2 \epsilon^2
\end{array}
$$
which completes the proof of the lemma.
\end{proof}

We finish the proof of the Theorem~\ref{relp} by noticing that the
measure $\mu_v$ satisfies all the conditions in Lemma~\ref{Rp},
and also that $\mu_v$ is supported on $\T^{p}\setminus \{0\}$.
This implies that
$$
(\sqrt{p+25} + 3)^2 \epsilon^2 \geq 1,
$$
which is equivalent to
$$
\epsilon \geq \frac{1}{\sqrt{p+25}+3}.
$$
Therefore the first inequality is not satisfied if $\epsilon \leq
1/l(p)$. This proves that, if
the representation $(\pi, \mr{H})$, does not have $\Z^{p}$
invariant vectors, then for any $v$, there exists $g\in F\cup G$
such that $||\pi(g)v-v|| \geq \frac{||v||}{l(p)}$.
\end{proof}
\begin{co}
\label{abelian1}
Let $(\pi ,\mr{H})$ be a unitary representation of the group
$$
G= \SL_p(\Z)\ltimes \Z^{p}.
$$
Let $v\in \mr{H}$ be a
$(F \cup G,\epsilon)$ invariant vector.
Then for every $g$ in $\Z^{p}$
we have $||\pi(g)v - v || \leq 2l(p)\epsilon$.
\end{co}
\begin{proof}
Let us split the Hilbert space $\mr{H}$ as a direct sum of the
closed subspaces $\mr{H}_0$ and $\mr{H}_1$, where $\mr{H}_0$
contains all $\Z^{p}$ invariant vectors and $\mr{H}_1$ is the
orthogonal compliment of $\mr{H}_0$. We have that both $\mr{H}_0$
and $\mr{H}_1$ are closed under the action of the group $G$,
because $\Z^{p}$ is a normal subgroup of $G$. Lets us write
$v=v_0+v_1$, where $v_i \in \mr{H}_i$. Since there are no $\Z^{p}$
invariant vectors in $\mr{H}_1$, there
exists $h \in F \cup G$ such that %
$||\pi(h)v_1 - v_1|| \geq ||v_1||/l(p)$. But we have that
$$
||\pi(h)v - v||^2 = ||\pi(h)v_0 - v_0||^2 + ||\pi(h)v_1 - v_1||^2
    \leq \epsilon^2,
$$
therefore $||v_1|| \leq l(p) \epsilon$. For any $g\in \Z^{p}$, we
have
$$
||\pi(g)v - v||^2 = ||\pi(g)v_0 - v_0||^2 + ||\pi(g)v_1 - v_1||^2
   \leq 0 + 4 || v_1||^2 \leq 4(l(p)\epsilon)^2,
$$
therefore $||\pi(g)v - v || \leq 2l(p)\epsilon$.
\end{proof}
\end{section}

\begin{section}{Relative Kazhdan constant for
$(\SL_p\times \SL_q)\ltimes \Z^{pq}$}
\label{sec:relpq}

In this section we estimate the relative Kazhdan constant for the
maximal parabolic subgroup %
$\left(\SL_p(\Z)\times \SL_q(\Z) \right)\ltimes \Z^{pq}$ of
$\SL_{p+q}$, with respect to the group $\Z^{pq}$ considering the
set of elementary matrices in $\SL_p$ and $\SL_q$ together with
the basis vectors of $\Z^{pq}$ as a generating set. The proof of
Theorem~\ref{rel} is based on Theorem~\ref{relp}.

\begin{thq}
\label{rel}
Let $F_1$ and $F_2$ denote the sets of elementary matrices in
$\SL_p(\Z)$ and $\SL_q(\Z)$ respectively, and
$G$ denotes the set of the $pq$ standard basis elements of $\Z^{pq}$.
Let $(\pi,\mr{H})$ be
a unitary representation of %
$(\SL_p(\Z)\times \SL_q(\Z))\ltimes \Z^{pq}$,
containing a vector $v$ which is $(F_1\cup F_2\cup G, 1/k(p+q))$
invariant, where %\marginpar{check it!}
$$
k(n) = \sqrt{5n/2+60}+6 .
$$
Then $\mr{H}$ contains $\Z^{pq}$ invariant vector, provided that
$p,q\geq 2$.
\end{thq}

\begin{rem}
Let us consider the standard representation of $\SL_{p+q}$ on
$\cl^2(\Z^{p+q})$. The group $(\SL_p\times \SL_q)\ltimes \Z^{pq}$
is isomorphic to a maximal parabolic in $\SL_{p+q}$, thus we have
a representation of that group in $\cl^2(\Z^{p+q})$. The
representation decomposes as a sum of two representations, one is
isomorphic to $\cl^2(\Z^p)$, where $\SL_q$ and $\Z^{pq}$ act
trivially, the other is a representation on %
$\cl^2(\Z^{p+q} \setminus \Z^p)$, without $\Z^{pq}$ invariant
vectors. Using that representation and a suitable vector $v$, it
can be shown that the Kazhdan constant is at most $\sqrt{2/q}$.
This shows that $1/\tilde{k}(n)$ is not a relative Kazhdan
constant, where $\tilde{k}(n) =\sqrt{n}/2$.
\end{rem}
\iffalse
\begin{rem}
If we simply generalized the proof by Burger, we will get a bound
of the type $[pq(p+q)]^{-1/2}$. This is true because if we
consider the action of $\SL_p \times \SL_q$ on the projective
space $P\R^{pq-1}$, there exists probability measure $\mu$ such
that for every Borel set $B$ we have $|\mu(gB) - \mu(B)| \leq
1/(p+q)$.
\end{rem}
\fi
\begin{proof}
%The proof is generalization of the similar result for the group
%$\SL_2(\Z) \ltimes \Z^2$, which can be found in~\cite{Shalom1}.

Assume that $v\in\mr{H}$ is $(F_1 \cup F_2\cup G, \epsilon)$
invariant vector, and that $\mr{H}$ does not contain $\Z^{pq}$
invariant vector.
%Let us restrict the representation $\pi$ to the $\Z^{pq}$ and
Let $P$ be the % corresponding
projection valued measure on $\widehat{\Z^{pq}} = \T^{pq}$, coming
form the restriction of the representation $\pi$ to $\Z^{pq}$, and
let $\mu_v$ be the measure on $\T^{pq}$, defined by %
$\mu_v(B)=\la P(B)v,v \ra$. The measure $\mu_v$ is supported on
$\T^{pq}\setminus \{0\}$, because by assumption $\mr{H}$ does not
contain an $\Z^{pq}$ invariant vectors.
%which implies that $P(\{0\}) = 0$.

We can identify the torus $\T^{pq}$ with the product of $q$ tori
of dimension $p$. For an element $x \in \T^{pq}$ we will write
$x=(x_1,\dots,x_q)$, where each $x_i$ is in $\T^p$, we will also
write $x_i= (x_{i1},\dots,x_{ip})$, where all $x_{ij}$ are in
$\R/\Z$, which we identify with the interval $(-1/2,1/2]$.

The proofs  of the next two lemmas are similar to the ones in
section~\ref{relp} and we will omit their proofs.

\begin{lm}
Let $K_{ij} = \{x \mid 1/4 > |x_{ij}| \}$, then %
$\mu_v(K_{ij}) \geq 1-\epsilon^2/2$.
\end{lm}

%\begin{rem}
%Let $K_{i} = \{x \mid |x_{ij}| <1/4 \}$ then
%$\mu_v(K_i) \geq 1-p\epsilon^2/2$,
%because the compliment of the union of the
%sets $K_{ij}$ (when we vary $j$)
%is exactly $K_i$.
%\end{rem}

\begin{lm}
\label{invariance1}
For every Borel set $B \subset \T^{pq}$ and
every elementary matrix $g \in F_1 \cup F_2$, we have that
$$
|\mu_v(gB) - \mu_v(B)| \leq 2\epsilon\sqrt{\mu_v(B)}+ \epsilon^2.
$$
The action of $\SL_p(\Z) \times \SL_q(\Z)$ on $\T^{pq}$, comes
from the standard action on $\R^{pq}$, by the isomorphism
$\T^{pq}=\R^{pq}/\Z^{pq}$.
\end{lm}

We need a result similar to Lemma~\ref{Rp}, considering the action
of $\SL_p \times \SL_q$ on $\T^{pq}$.
\begin{lm}
\label{Rpq}
Let $\mu$ be a finitely additive probability measure on $\T^{pq}$
%(considered as a subset of $\R^{pq}$)
such that
\begin{itemize}
\item %
$\mu(|x_{ij}| \geq 1/4) \leq \epsilon^2/2$ for any $i$ and
$j$, %
\item %
$|\mu(gB) - \mu(B)| \leq 2\epsilon \sqrt{\mu(B)} + \epsilon^2$ %
for any Borel set $B$ and any elementary matrix $g$ in $\SL_p(\Z)$
or $\SL_q(\Z)$.
\end{itemize}
Then the measure of the origin is at least
$$
\mu(\{x \mid x=0\}) \geq
1 -  (\sqrt{3p+2q+60}+6)^2\epsilon^2.
$$
\end{lm}
\begin{proof}
For a point $x\in \T^{pq}$ we write $x=(x_{1}\dots,x_{q})$, where
$x_i \in \T^p$ and $x_i =(x_{i1},\dots,x_{ip})$ and $x_{ij}\in
(1/2,1/2]$. Let us define the Borel sets
$$
\begin{array}{l}
B_i = \{ x \mid x_k =0 \mbox{ for } k\leq i\}, \mbox{ and}\\
C_i = \{ x \mid x_1 = x_i \not =0, x_k = 0 \mbox{ for } 1<k<i\}.
\end{array}
$$
The elementary matrix $g_{1i} \in \SL_q$ sends %
$B_{i-1} \setminus B_i$ into $C_i$ for any $i\geq 2$. Therefore, we
have
$$
\mu(B_{i-1} \setminus B_i ) \leq
    \mu(C_i) + \epsilon^2 + 2\epsilon\sqrt{\mu(C_i)}.
$$
Let us notice that the sets $C_i$ are disjoint and lies in the
compliment of $B_1$. Therefore
$$
\mu(B_1 \setminus B_q ) \leq
     \mu(\T^{pq}\setminus B_1) + (q-1)\epsilon^2 +
     2\epsilon\sqrt{(q-1) \mu(\T^{pq}\setminus B_1)}.
$$
Using lemma~\ref{Rp} (by considering the measure $\tilde \mu$
on $\T^p$, defined as follows: %
$\tilde \mu(K) = \mu(\{x | \mid x_1 \in K \})$) we have
$$
\mu(\T^{pq}\setminus B_1)\leq
    (\sqrt{p+25}+3)^2\epsilon^2
$$
Finally we have
$$
\begin{array}{r@{\,}l}
\!\!\!\mu(\T^{pq} \setminus B_q) \!\leq & \displaystyle
    2\mu(\T^{pq}\setminus B_1) +
        2\epsilon\sqrt{(q-1) \mu(\T^{pq}\setminus B_1)} +
        (q-1)\epsilon^2 \leq \\
\leq & \displaystyle
    (2p + 66 + 12\sqrt{p})\epsilon^2 +
    2\epsilon^2\sqrt{(q-1)(p+33+6\sqrt{p})} +
    (q-1)^2\epsilon^2 \!\!\leq \\
\leq & \displaystyle
    (3p + 2q + 97 + 18\sqrt{p})\epsilon^2 \leq
    (\sqrt{3p+2q+60}+6)^2\epsilon^2.
\end{array}
$$
That completes the proof of the lemma since the set
$B_q$ contains only the origin.
\end{proof}

We finish the proof of the theorem considering without loss of
generality that $q\geq p \geq 2$. The measure $\mu_v$ satisfies
all the conditions in the lemma~\ref{Rpq}, and also $\mu_v$ is
supported on $\T^{pq}\setminus \{0\}$, because $\mr{H}$ does not
have invariant vectors. This implies that
$$
(\sqrt{3p+2q+60}+6)^2\epsilon^2 \geq 1,
$$
which is equivalent to
$$
\epsilon \geq \frac{1}{\sqrt{3p+2q +60}+6}.
$$
This inequality is not satisfied if $\epsilon \leq 1/k(p+q)$.

This proves that if the representation $(\pi, \mr{H})$, does not
have $\Z^{pq}$ invariant vectors, then for any $v$, there exists
$g\in F_1\cup F_2 \cup G$ such that %
$||\pi(g)v-v|| \geq \frac{||v||}{k(p+q)}$.
\end{proof}
\iffalse
\begin{rem}
This proof works only when $p,q \geq 2$. If one of them is $1$,
we can use theorem~\ref{relp} which gives a slightly better bound.
\end{rem}
\fi

\begin{co}
\label{abelian}
Let $(\pi ,\mr{H})$ be a unitary representation of the group
$$
G= (\SL_p(\Z)\times \SL_q(\Z))\ltimes \Z^{pq}.
$$
Let $v\in \mr{H}$ be a
$(F_1\cup F_2\cup G,\epsilon)$ invariant vector.
Then for every $g$ in $\Z^{pq}$
we have $||\pi(g)v - v || \leq 2k(p+q)\epsilon$.
\end{co}
%\iffalse
\begin{proof}
The proof of this corollary is similar to the one of corollary~\ref{abelian1}
\end{proof}
%\fi
\end{section}

\begin{section}{Vectors systems in $\Z^k$}
\label{sec:vec}

Let $v_1,\dots,v_n$ be vectors in $\Z^k$, which generate the whole
group $\Z^k$, we will call $V=\{v_1,\dots,v_n\}$ a complete system
of vectors in $\Z^k$. We can also consider $V$ as a left
invertible $k\times n$ matrix with integer coefficients by letting
$V=(v_1,\dots,v_n)^t$.

We can define an elementary transformation $E_{i,j,a}$ on a
complete vector system $V$, which preserves all vectors except
$v_j$ and sends $v_j$ to $v_j' = v_j + a v_i$. It is clear the we
obtain a new complete vector system after this operation.

It is well known (see~\cite{Carter,Carter2,Vas}) that if %
$n \geq k+2$,%
\footnote{The condition $n\geq k+2$, comes from the fact
that the ring $\Z$ has stable range equal to $2$, see~\cite{Vas}
for details. Using the fact that $\SL_3(\Z)$ is boundedly generated
by the elementary matrices it is possible to extend this result
to all $n \geq k$ except $n=2$ and $k=1$ or $k=2$.} %
using approximately $2kn$ elementary operations we can transform
any vector system to the canonical vector system $U$, which
contains only standard basis vectors at the first $k$ places and
the zero vectors in the other places.

In this section we will show that using a few `generalized
elementary transformations' we can also transform any complete
vector system $V$ to the canonical one $U$.

Let us partition the set of indices $\{1,\dots,n\}$ into two
disjoint parts $I$ and $J$. For any $|I|\times |J|$ matrix
$\alpha$ we define a generalized elementary transforation
$E_{I,J,\alpha}$ as follows: For any vector system $V=\{v_i\}$, we
define a new vector system $V'=\{v_i'\}$ as follows:
\begin{list}{}{}%
\item %
$v_i' = v_i$ for all $i\in I$; %
\item %
$v_j' = v_j + \sum_{i\in I} \alpha_{ij}v_i$ for all $j\in J$.
\end{list}

If we consider $V$ as a $k \times n$ matrix with integer
coefficients, the generalized elementary operation
$E_{I,J,\alpha}$ corresponds to left multiplication with the
matrix $A$, obtained from %
$\left(\begin{array}{cc}
I & \alpha \\ 0 & I
\end{array}\right)$, %
by rearranging the rows and the columns.

\iffalse It is well known that if $n \geq k+2$, using
approximately $2kn$ elementary operations we can transform any
vector system to the standard one, which contains only the
standard basis vectors and the zero ones. The next theorem shows
that we can do the same by using just a few generalized elementary
operations.
\fi

\begin{thq}
\label{vectors}
If $n\geq 3k$, then any complete system $V$of $n$ vectors in $\Z^k$
can be transformed by using at most $4$ generalized elementary operations,
to the `standard' system of vectors $U$, where the first $k$ vectors in $U$
are the standard basis vectors of $\Z^k$ (in the same order)
and all other vectors are zero.
\end{thq}
\begin{rem}
This result for $k=1$ is well known and it is used in the
induction step of the proof that $\SL_n(\Z)$ is bounded generated
by elementary matrices. In fact for $k=1$, three operations are
enough. For $k>1$ using $3$ operations we know how to transform
the system $V$, into a system $U'$, which contains $k$ vectors
from the standard basis of $\Z^k$ and $n-k$ zero vectors, but we
do not know how to control the positions of the nonzero vectors.
\end{rem}
\iffalse
\begin{rem}
The condition $n\geq 3k$ is too strong and can be replaced by
$n\geq 2k+1$, but in that case we need $5$ operations. This is
true because for any `prime' subgroup $B$ the quotient $\Z^k/B$ is
cyclic and can be generated by 1 element. So using an additional
generalized linear transformation (after the first in the proof),
we can modify $1$ vector so that some $k+1$ vectors generate the
whole group $\Z^k$%
\footnote{This observation was made by Tal Poznansky, it is used
for obtaining the bound in Theorem~\ref{main:strong}.}. %
Also if we replace $\Z$ with some field than $3$ generalized
elementary transformations are enough provided that $n\geq 2k$.
\iffalse

Its is interesting wether this condition can be replaced by %
$n \geq k+C$, for some fixed constant $C$. Such result will
improve the Kazhdan constant for $\SL_n(\Z)$ by approximately a
factor of approximately $3$ (if the number of transformations
stays the same).
\fi
\end{rem}
\fi

\begin{proof}
Let us first recall how to transform such system in the case
$k=1$. By one elementary operation we can make one of the vectors
a sufficiently big prime number, after another operation we can
put $1$ at the first place, and using the final operation we can
make all other vectors equal to $0$.

In order to generalize this construction for $k\geq 2$ we need to
define the analog of the prime number.
\begin{df}
We call a finite index subgroup $B$ of $\Z^k$
%generated by linearly independent vectors $p_1,\dots,p_k$
a `prime' subgroup if
the quotient $\Z^k/B $ is isomorphic to
$$
\Z/\pi_1\Z \times \dots \times \Z/\pi_k\Z,
$$
where $\pi_i$ are pairwise different prime numbers. Any $k$
vectors which generate a `prime' subgroup are called a `prime'
system.
\end{df}
\begin{rem}
The vectors $w_1=(\pi_1,*,*,\dots,*)$, $w_2=(0,\pi_2,*,\dots,*)$,
%$w_3=(0,0,\pi_3,\dots ,*)$,
\dots, $w_k=(0,0,0,\dots,\pi_k)$, where $\pi_i$ are distinct primes
generate a prime subgroup of $\Z^k$.
\end{rem}

\begin{lm}
Let $V$ be a system of vectors in $\Z^k$, using one generalized
elementary operation we can transform $V$ into a system $V'$,
where some $k$ vectors form a `prime' system, i.e., they generate
a `prime' subgroup of $\Z^k$.
\end{lm}
\begin{proof}
Assume that the last $k$ vectors are linearly independent.
\begin{rem}
Using several elementary transformations (which modify only the
first $k$ vectors), we can transform any complete system $V$
into system $V'$ such that $v_1' = (\pi_1 ,*,*, \dots,*)$,
$v_2'= (0,\pi_2,*,\dots,*)$, \dots, $v_k'= (0,0,0,\dots,\pi_k)$,
where $\pi_i$ are sufficiently large distinct prime numbers
--- larger then the determinant of the matrix formed by the
coefficients of the last $k$ vectors.
\end{rem}
\begin{proof}
The proof is by induction on $k$ --- in the base case $k=0$, there
is nothing to prove. Suppose that the vectors $v_1,\dots,v_{k-1}$
have the desired form. Let us consider the set of vectors
$$
P=\{v_k + \sum_{i\not=k} \alpha_i v_i | \alpha_i \in \Z\}.
$$
Since the vectors $\{v_i\}$ form a complete vector system and the
last $k$ vectors are linearly independent, the set $P$ contains
all vectors of the form $(0,\dots,0,a+\lambda .d)$, for all
$\lambda \in \Z$, for some relatively prime integers $a$ and $d$.
Here we use that $\pi_i$, for $i< k$, are sufficiently big prime
numbers therefore the standard basis vectors $e_i$ for $i< k$ lie
in the subgroup generated by the vectors %
$ v_1,\dots, v_{k-1}, v_{n-k+1},\dots, v_n$ .

Using Dirichlet's theorem about primes in the arithmetic
progressions, it follows that $P$ contains the vector of the form
$ (0,0,0,\dots,\pi_k)$, which completes the induction step.
\end{proof}

By the above remark using elementary transformations which modify
only the first $k$ vectors we can make these vectors a `prime'
system. Doing all these elementary transformations corresponds to
multiplying from the left (the matrix of the vector system $V$)
% with coordinates of the vectors $v_i$)
with matrix $A \in \SL_n(\Z)$ of the form
$$
\left( \begin{array}{cc}
* & * \\ 0 & I
\end{array}\right),
$$
where the blocks are of sizes $k$ and $n-k$. Any such matrix can
be written uniquely as $A = B C$, where $B,C \in \SL_n(\Z)$ and
$$
B = \left( \begin{array}{cc}
* & 0 \\ 0 & I
\end{array}\right)
\quad
C = \left( \begin{array}{cc}
I & * \\ 0 & I
\end{array}\right).
$$
If we apply the generalized elementary transformation
corresponding to the multiplication by matrix $C$, we obtain
vector system such that the subgroup generated by the first $k$
vectors, coincides with the subgroup generated by $v_i'$-es.
Because the upper left corner of $B$ is in $\SL_k(\Z)$, and the
multiplication by $B$ does not change the subgroup generated by
the first $k$ vectors, this subgroup is `prime', which
finishes the proof of the lemma.
\end{proof}

Notice that if $B$ is a `prime' subgroup in $\Z^k$, then any
strictly increasing sequence of subgroups between $B$ and $\Z^k$
has at most $k$ terms. This implies that if a complete vector
system contains a `prime' subsystem, then there exist at most $2k$
vectors which generate the whole group $\Z^k$.

For such system by applying one generalized elementary operation
we can generate $k$ vectors that form a standard basis of $\Z^k$. This
is true, because we have $2k$ vectors, which generate the whole
$\Z^k$ and putting them in the set $J$, we can transform the other
vectors to any vectors in $\Z^k$. Since $n\geq 3k$, we have at
least $k$ vectors to modify and we can make these vectors equal to
the standard basis vectors $e_i$. Moreover, if any of these vectors
is among the first $k$ we can make it equal to the corresponding
vector in the standard basis.

Finally, we need one more transformation in order to make the
first $k$ vectors equal to `standard' basis vectors of $\Z^k$ and
with one final generalized elementary operation we can make all
the other vectors $0$-es.
\end{proof}
\begin{rem}
The condition $n\geq 3k$ is too strong and can be replaced by
$n\geq 2k+1$, but in that case we need $5$ operations. This is
true because for any `prime' subgroup $B$ the quotient $\Z^k/B$ is
cyclic and can be generated by 1 element. So using an additional
generalized linear transformation (after the first in the proof),
we can modify $1$ vector so that some $k+1$ vectors generate the
whole group $\Z^k$%
\footnote{This observation was made by Tal Poznansky, it is used
for obtaining the bound in Theorem~\ref{main:strong}.}. %
Also if we replace $\Z$ with some field then $3$ generalized
elementary transformations are enough, provided that $n\geq 2k$.

It is interesting whether this condition can be replaced by %
$n \geq k+C$, for some fixed constant $C$. Such a result will
improve the Kazhdan constant for $\SL_n(\Z)$ by approximately a
factor of $3$ (if the number of transformations
stays the same).
\end{rem}

\begin{co}
\label{decompose}
If $n\geq 3k$, then
any matrix $g\in\SL_n(\Z)$ can be written as a product of 6 matrices:
$$
g=g_1g_2g_3g_4g^*g_5,
$$
where $g^*$ lies in the copy of $\SL_{n-k}(\Z)$ embedded in the
lower right corner. Also any matrix $g_i$ can be obtained from a
matrix of the type %
$\left(\begin{array}{cc} %
I & * \\ 0 & I %
\end{array} \right)$, %
by rearranging the rows and columns
(%the sizes and
the position of the blocks depend on the matrix $g_i$).
\end{co}
\begin{proof}
Let us consider the first $k$ entries of each row of $g$. They
form a complete system of $n$ vectors in $k$ dimensional space,
because $g$ is an invertible matrix.

Every generalized elementary transformation on these vectors
corresponds to multiplying the $k\times n$ matrix of their
coordinates from the left by a matrix similar to %
$\left(\begin{array}{cc} %
I & * \\ 0 & I %
\end{array} \right)$. %
By Theorem~\ref{vectors} after 4 such multiplications we can
transform this matrix to %
$\left(\begin{array}{c} %
I \\ 0 %
\end{array} \right)$.

Therefore, by multiplying $g$ from the left with these matrices we
can reduce it to a matrix of type %
$\left(\begin{array}{cc} %
I & * \\ 0 & * %
\end{array} \right)$. %
Finally by one multiplication from the right we can transform this
matrix to %
$g^*=\left(\begin{array}{cc} %
I & 0 \\ 0 & * %
\end{array} \right)$, %
which lies in $\SL_{n-k}$. If we `reverse' this process we
obtained the desired decomposition of the matrix $g$.
\end{proof}
\begin{co}
If $n\geq 3$, then any matrix $g\in\SL_n(\Z)$ can be written as a
product of at most $60 + 13 \ln n$ matrices, each of which can be
obtained from a matrix of the type %
$\left(\begin{array}{cc} %
I & * \\ 0 & I %
\end{array} \right)$, %
by rearranging the rows and columns.
\end{co}
\end{section}

\begin{section}{Kazhdan constants for $\SL_n(\Z)$}
\label{sec:proof}

Using the fact that $\SL_n(\Z)$ (for $n\geq 3$) is bounded
generated by the elementary matrices, and using an analog of
Corollary~\ref{abelian} for $\SL_2(\Z) \ltimes \Z^2$, it can be
shown (see~\cite{Shalom1}), that if $(\pi,\mr{H})$ is an unitary
representation of $\SL_n(\Z)$, and $v$ is an $\epsilon$ invariant
vector with respect to all elementary matrices, then for any %
$g \in \SL_n(\Z)$ we have that %
$||\pi(g)(v) - v || \leq 22 f(n)\epsilon$, where %
$f(n) = 3(n^2 - n)/2 + 51$, which is the number of elementary
matrices (with any integer off the diagonal) needed to express any
element in $\SL_n(\Z)$, see~\cite{Carter2}. From here it easily
follows that the Kazhdan constant for $\SL_n(\Z)$ with respect to
the elementary matrices is at least $1/22f(n)$. Our goal is to
improve the upper bound  $22f(n)$ and obtain a better Kazhdan
constant.

\begin{df}
Let $h(n)$ is the smallest number such that for any unitary
representation and any positive number $\epsilon$, the condition
$||\pi(g)(v)-v||<\epsilon$ for any elementary matrix $g$ in $E_n$,
implies that $||\pi(g)(v)-v|| \leq h(n) \epsilon$ for any %
$g\in \SL_n(\Z)$.
\end{df}

\begin{lm}
\label{bound} If $n\geq 3i$ and $i\geq 2$, then the function
$h(n)$ satisfies the inequality
%\marginpar{check it }
$$
h(n) \leq h(n-i) + 10k(n) \leq h(n-i) + \sqrt{250n + 6000} + 60.
$$
Here $k(n)$ is the function defined in Theorem~\ref{rel}.
\end{lm}
\begin{proof}
Let $(\pi,\mr{H})$ be a unitary representation of $\SL_n(\Z)$ and
$v\in \mr{H}$ be a unit vector such that %
$||\pi(g)v-v|| \leq \epsilon$ for any elementary matrix $g$.

%For any matrix $A$ be a matrix of the form
%$\left(\begin{array}{cc} I & * \\ 0 & I \end{array} \right)$.
The set of all matrices of the form %
$\left(\begin{array}{cc} %
* & * \\ 0 & * %
\end{array} \right)$ %
is a subgroup of $\SL_n(\Z)$ isomorphic to %
$(\SL_p(\Z) \times \SL_q(\Z)) \ltimes \Z^{pq}$. If we restrict the
representation $\pi$ to this subgroup, we can apply
corollary~\ref{abelian} and obtain that %
$||Av-v|| \leq 2k(p+q)\epsilon \leq 2k(n)$ for any matrix $A$ of
the form %
$\left(\begin{array}{cc} %
I & * \\ 0 & I %
\end{array} \right)$.

Let $g$ be a matrix in $\SL_n(\Z)$ by lemma~\ref{decompose} we can
write $g$ as a product of $6$ matrices. By the above argument,
five of these matrices move the element $v$ by less than
$2k(p+q)\epsilon\leq 2k(n)\epsilon$. The sixth matrix lies in a
copy of the group $\SL_{n-k}(\Z)$ and if we restrict the
representation $\pi$ to that subgroup we can see that it moves the
vector $v$ by less than $h(n-i)\epsilon$. This implies that $g$
moves $v$ by less than $(h(n-i) + 10k(n))\epsilon$, which proves
the lemma.
\end{proof}

Before completing the proof of Theorem A, we need a lemma about
functions which satisfy an inequality like the one in
lemma~\ref{bound}.
\begin{lm}
\label{recursion} Let $a,b,c$ be positive real numbers,
$\lambda < 1$ and let $f:\N \to \R$ be a %non decreasing
function.
If the function $f$ satisfies the inequality
$$
f(n) \leq f(i) + \sqrt{an+b} + c,
$$
for any $i \geq \lambda^2 n$ and any %
$n \geq n_0 > 1/(1-\lambda^2)$ then
$$
f(n) \leq %
A(\sqrt{n} - \lambda \sqrt{\widetilde{n_0}}) - %
c\left(\log_{\lambda^2} \frac{n}{\widetilde{n_0}}+1\right)  + %
\frac{B}{\sqrt{\widetilde{n_0}}}+ f(n_0)
$$
where $A$, $B$ and $\widetilde{n_0}$ are given by:
$$
A=\frac{\sqrt{a}}{1-\lambda} \quad %
B= \frac{b + a/(1-\lambda^2)}{(1-\lambda)\sqrt{a}} \quad %
\widetilde{n_0} = n_0 - \frac{1}{1-\lambda^2}.
$$
\end{lm}
\begin{proof}
Let us define recursively the sequences $x_k$ and $y_k$ as
follows: $x_0=y_0=n$ and $x_{i+1} = \lceil \lambda^2 x_i \rceil$,
$y_{i+1} = \lambda^2 y_i$. Here $\lceil x \rceil$ denotes the
smallest integer greater than $x$. By induction it follows that
$$
y_i < x_i \leq %
y_i + \frac{1-\lambda^{2i}}{1-\lambda^2} %
< y_i + \frac{1}{1-\lambda^2}.
$$
Therefore, for %
$s =\lceil -\log_{\lambda^2} (n/\widetilde{n_0}) \rceil$ %
we have %
$x_s \leq y_s + 1/(1-\lambda^2) \leq %
\widetilde{n_0} + 1/(1-\lambda^2) = n_0. $ %
Using the functional inequality we have
$$
f(x_{i+1}) \leq f(x_i) + \sqrt{ax_{i+1} + b} +c.
$$
Adding all these inequalities for differen $i$'es we obtain
$$
f(n) \leq f(x_s) + \sum_{i=0}^{s-1} (\sqrt{ax_i + b} +c) \leq %
f(n_0) + \sum_{i=0}^{s-1} (\sqrt{ax_i + b} +c)
$$
Using the inequality between $x_i$ and $y_i$ we have
$$
\begin{array}{r@{\,\,}l}
f(n) \leq &  \displaystyle %
    f(n_0) + c s + %
        \sum_{i=0}^s\sqrt{an\lambda^{2i} + a/(1-\lambda^2) +b} = %
\\ %
= & \displaystyle %
    f(n_0) + c s + \sqrt{an} \sum_{i=0}^s \lambda^i %
        \sqrt{1 + \lambda^{-2i}\frac{a/(1-\lambda^2) +b}{an}} \leq %
\\ %
\leq &  \displaystyle %
    f(n_0) + cs + \sqrt{an} \sum \lambda^i + %
        \sum \frac{a/(1-\lambda^2) +b}{2\sqrt{an}}\lambda^{-i} \leq %
\\
& \quad \mbox{because } %
    \sqrt{1+x} \leq 1 + x/2 \mbox{ for every $x$} %
\\
\leq &  \displaystyle %
    f(n_0) + cs  + %
        \frac{\sqrt{an}(1-\lambda^s)%
        \rule[15pt]{0pt}{0pt}}{1-\lambda} + %
        \frac{a/(1-\lambda^2) +b}{\sqrt{an}} %
        \frac{\lambda^{-s}}{1-\lambda} \leq %
\\
\leq & \displaystyle %
    f(n_0) + cs + A\sqrt{n}(1-\lambda^s)  + %
        B \lambda^{-s} /\sqrt{n} \leq %
\\
\leq & \displaystyle %
    A(\sqrt{n} - \lambda\sqrt{n_0}) + %
        c\left(\log_{\lambda^2} \frac{n}{\widetilde{n_0}}+1\right) %
        + \frac{B}{\sqrt{\widetilde{n_0}}}+ f(n_0),
\end{array}
$$
For the last inequality we used %
$s \leq \log_{\lambda^2} \frac{n}{\widetilde{n_0}}+1$ and %
$\lambda^{-2s} \leq n/\widetilde{n_0}$.
\end{proof}

Applying the previous lemma to the function $h(n)$ we obtain
\begin{thq}
The function $h(n)$ satisfies the inequality
%\marginpar{Check the constant!}
$$
h(n) < 90 \sqrt{n} + 4000.
$$
\end{thq}
\begin{proof}
By lemma~\ref{bound} we have that the function $h(n)$ satisfies the
inequality with $a=250$, $b=6000$, $c=60$ and %
$\lambda = \sqrt{2/3}$. Putting  these constants and $n_0 = 7$ in
the lemma~\ref{recursion} gives
$$
A = \frac{\sqrt{250}}{1-\sqrt{2/3}} = 15\sqrt{10}+10\sqrt{15} \quad % 86.16
B = \frac{6750}{\sqrt{150}(1-\sqrt{2/3}))} = 675(2+\sqrt{6}), %3003.41
$$
and $\widetilde{n_0} \geq 1$, which implies the inequality
$$
h(n) \leq  %
    (15\sqrt{10}+10\sqrt{15})(\sqrt{n} - \sqrt{14/3%\frac{14}{3}
    })
    +  60(\log_{3/2%\frac{3}{2}
    } \frac{n}{7} + 1)  + 675(2+\sqrt{6})/\sqrt{7}  %
    + h(7).
$$
By the Shalom result we have $h(n) \leq 33n^2-11n+1152$, i.e.,
$h(7) \leq 2692$. Finally we have
$$
h(n) < (15\sqrt{10}+10\sqrt{15})\sqrt{n} +  60\log_{3/2} n + 3900 < %
    90\sqrt{n} + 4000.
$$
\end{proof}

Now we prove Theorem~A.
{\renewcommand{\thelm}{A}\addtocounter{lm}{-1}
\begin{thq}
The Kazhdan constant for $\SL_n(\Z)$ and $\SL_n(\F_p)$ with
respect to the elementary matrices is
$$
\mr{K}(\SL_n(\Z),E_n) \geq (64\sqrt{n}+2850)^{-1}.
$$
\end{thq}
}
\begin{proof}
It is well known fact %(see for example ..)
that if a representation $(\pi, \mr{H})$ of a group $G$ contains a
unit vector $v\in \mr{H}$ such that $|| \pi(g)v - v|| < \sqrt{2}$
for any $g\in G$ then $\mr{H}$ contains a $G$-invariant vector.
Applying this observation gives that
$$
\mr{K}(\SL_n(\Z),E_n) \geq \sqrt{2}/h(n) %
    \geq (50\sqrt{n}+2850)^{-1}.
$$
\end{proof}

\begin{rem}
More detailed consideration, using the exact size of the blocks of
matrices $g_i$ in Corollary~\ref{decompose}, and using the
stronger version of Theorem~\ref{vectors} for $n\geq 2k+1$, gives
that
$$
h(n) \leq %
   \sqrt{2}(5\sqrt{5}+1)(\sqrt{2}+1)\sqrt{n} + 22\log_2 n + 350 < %
   \sqrt{2}( 42 \sqrt{n}+ 860),
$$
which implies the bound of the Kazhdan constant in
Theorem~\ref{main}
$$
\mr{K}(\SL_n(\Z),E_n) \geq (33\sqrt{n} + 317)^{-1}.
$$
Similarly using the version of Theorem~\ref{vectors}, for vector
systems over a field $\F_p$ we have
$$
h(n) \leq 8\sqrt{3}(\sqrt{2}+1)\sqrt{n} + 8\log_2 n /3 + 100 <
\sqrt{2}( 24 \sqrt{n}+ 100),
$$
which proves Theorem~\ref{main:fp}.
\end{rem}
\end{section}

\begin{section}{Generalizations to other groups}

In this section we show how Theorem~\ref{main} can be generalized
to the groups $\SL_n(R)$ for several classes of rings $R$. We will
only sketch the proofs of the necessary lemmas. In order to do so
we need to generalize Theorem~\ref{rel} and Lemma~\ref{decompose}.

The first step is the proof analogous to the proof of the
lemma~\ref{R2} for the
ring $\Z[t_1,\dots,t_s]$. Let $F$ denote the set of elementary
matrices in $\SL_2$, with $\pm 1$ and $\pm t_i$ off the diagonal.
\begin{lm}
\label{R2s} Let $\mu$ be a finitely additive measure on the dual
of $\Z[t_1,\dots,t_s]^2$, i.e. %
$\widehat{\Z[t_1,\dots,t_s]}^2= %
(\R/\Z[[t_1^{-1},\dots,t_s^{-1}]])^2$ %
such that
\begin{itemize}%
\item %
$\mu(|x_0| \geq 1/4) \leq \epsilon^2/2$ and %
$\mu(|y_0| \geq 1/4) \leq \epsilon^2/2$. %
Here $x_0$ denotes the constant term of the series $x$;

\item %
$|\mu(gB) - \mu(B)| \leq 2\epsilon \sqrt{\mu(B)} + \epsilon^2$ %
for any Borel set $B$ and any elementary matrix %
$g\in F \subset \SL_2(\Z[t_1,\dots,t_s])$.
\end{itemize}
Then if $\epsilon < 1/12$ we have
%\marginpar{check!!!}
$$
%???
\mu(\widehat{\Z[t_1,\dots,t_s]}^2\setminus \{(0,0)\}) \leq %
10.6^k \epsilon^2
$$
\end{lm}
\begin{proof}
The proof is by induction using lemma~\ref{R2} as the base case.
The proof of the induction step uses the description of
$\widehat{R[t]}$ in term of $\widehat{R}$ and is based on the
proof of lemma 3.3 from~\cite{Shalom1}.
%induction \marginpar{sketch it}
\end{proof}

Using this lemma we can generalize Theorems~\ref{relp}
and~\ref{rel} for the ring $\Z[t_1,\dots,t_s]$, and therefore for
any finitely generated ring, by replacing the functions $l(n)$ and
$k(n)$ with
%\marginpar{check!!!}
$$
%???
l_s(n) = \min\{\sqrt{3n+21.6^s},12\} \quad %
k_s(n)=\min\{\sqrt{6n+48.6^s},12\}
$$
This can be further generalized to rings which contain a
finitely generated dense sub rings -- like $\C$ or
$\Z[[t_1,\dots,t_s]]$.
\medskip

Theorem~\ref{vectors} can be generalized to many different classes
of rings replacing $\Z$. In the proof of Theorem~\ref{vectors}, we
used the fact that for any ideal $I \lhd \Z$ and any $x$ in the
ring $\Z$, such that $x \Z  + I = \Z$, there are infinitely many
elements $y \in x + I$, such that the ring $Z/yZ$ has a unique
maximal ideal. Therefore, for any commutative ring $R$ which has
the above property, Theorem~\ref{vectors} holds, and any vector
system in $R^k$ consisting of more than $3n$ vectors can be
transformed to the standard one using at most $4$ generalized
elementary transformations. An example of a ring satisfying this
condition is $\Z[[t]]$.

Suppose that the ring $R$ satisfies the following condition: there
are invertible elements in the coset $x + I$ for any element $x\in
R$ and ideal $I\lhd R$, such that $xR + I = R$, in particular if
$R$ is a local ring or a filed. Then we can transform any vector
system in $R^k$ to the standard one using $3$ generalized
elementary transformations of a fixed type, provided that $n \geq
2k$. A nontrivial example of a ring satisfying the above condition
is $\Hom(S^1,\C)$ with
pointwise operations%
\footnote{These conditions imply the stable range of the ring $R$
is at most $2$.}.

Using this remark we can show that if the ring $R$ satisfies one
of the above conditions and contains a dense sub-ring $S$
generated by $\alpha_i$ for $i=1,\dots,s$. Then the group
$\SL_n(R)$ has property $T$ and the Kazhdan constant is
%\marginpar{do it!!!}
$$
\mr{K}(\SL_n(R),E_n(R)) \geq (50\sqrt{n} + (10N+...) 6^s + 300)^{-1},
$$
provided that $\SL_3(R)$ is boundedly generated by the elementary
matrices, and every element $g \in \SL_3(R)$ can be written as a
product of $N$ elementary matrices. In particular we have the
following corollaries:

\begin{co}
For any be compact ring $R$ such that there exist $d$ elements which generate a
dense sub-ring, the groups $\SL_n(R)$, and
$\SL_n(R[[t_1,\dots,t_s]])$ have property $T$ and the Kazhdan
constant with respect to the set of elementary matrices is
$O(n^{1/2})$.
\end{co}
\begin{thq}
The loop group $\mr{L}(\SL_n(\C) = \SL_n(\mr{L}(\C))$ of
$\SL_n(\C)$, has property $T$ for $n\geq 3$ and the Kazhdan
constant with respect to the set of trivial loops $E_n$ of
elementary matrices with $\pm 1$ off the diagonal, is at least
%\marginpar{do it!!!}
$$
\mr{K}(\mr{L}(\SL_n(\C)),E_n) \geq %
 [50\sqrt{n} + B ]^{-1},
$$
where $B$ is a constant. Note that this is not a locally compact
Lie group and the set $E_n$ generates a finite dimensional
subgroup.
\end{thq}
\begin{proof}
The proof is based on the fact that the ring $\Hom(S^1,\C)$ contains a
dense sub-ring generated by $4$ elements -- $1$ and $\alpha_i$.
More over the elements $\alpha_i$ can be chosen in any
neighborhood of $0$, which allows us not to include in the generating set of
the group the
elementary matrices with $\pm \alpha_i$ off the diagonal.
\end{proof}

{\renewcommand{\thelm}{C}\addtocounter{lm}{-1}
\begin{thq}
Let $\mr{O}$ be the ring of integers in a number field
$\mathbb{K}$, with discriminant $\Delta$, which is generated as a
ring by $1$ and the elements $\alpha_i$ for $i=1,\dots,s$. The
Kazhdan constant for $\SL_n(\mr{O})$ with respect to the set
$E_n(\mr{O}$ of  elementary matrices with $\pm 1$ and %
$\pm \alpha_i$ is
$$
\mr{K}(\SL_n(\mr{O}),E_n(\mr{O})) \geq %
[50\sqrt{n} + (B + C\Delta)6^m]^{-1}.
$$
where $B$ and $C$ are universal constants and $\Delta$ is the
number of different prime divisors of the discriminant of
$\mathbb{K}:\Q$.
\end{thq}
}
\begin{proof}
Here we used the result by Carter and Keller~\cite{Carter} that
every element in $\SL_3(\mr{O})$ can be written as a product of
%\marginpar{check!!!}
$60 + \Delta$ elementary matrices.
\end{proof}

\end{section}

%\newpage

\bibliographystyle{plain}
\bibliography{SL}

\noindent
Martin Kassabov:\\
Department of Mathematics \\
University of Alberta\\
632 Central Academic Building\\
Edmonton, Alberta, T6G 2G1\\
Canada\\
E-mail: kassabov@aya.yale.edu

\end{document}